\documentstyle[12pt]{article}

 2
\def\BI{{\rm 1\!l}}
\def\pois#1#2{\{#1,#2\}}

\def\Eq#1{{\begin{equation} #1 \end{equation}}}

\def\gone{\matrix{{}\cr g \cr {}^1}}
\def\gtwo{\matrix{{}\cr g \cr {}^2}}
\def\hone{\matrix{{}\cr h \cr {}^1}}
\def\htwo{\matrix{{}\cr h \cr {}^2}}
\def\gbone{\matrix{{}\cr \bar g \cr {}^1}}
\def\gbtwo{\matrix{{}\cr \bar g \cr {}^2}}
\def\psone{\matrix{{}\cr \psi \cr {}^1}}
\def\pstwo{\matrix{{}\cr \psi \cr {}^2}}
\def\xone{\matrix{{}\cr x \cr {}^1}}
\def\xtwo{\matrix{{}\cr x \cr {}^2}}
\def\xtone{\matrix{{}\cr \tilde x \cr {}^1}}
\def\xttwo{\matrix{{}\cr \tilde x \cr {}^2}}
\def\ptone{\matrix{{}\cr \tilde p \cr {}^1}}
\def\pttwo{\matrix{{}\cr \tilde p \cr {}^2}}
\def\fone{\matrix{{}\cr f \cr {}^1}}

\begin{document}
\begin{flushright}
UAHEP 951\\
January 1995\\
\end{flushright}

\centerline{ \LARGE Lorentz transformations as Lie-Poisson symmetries}
\vskip 2cm

\centerline{ {\sc A. Simoni*, A. Stern** and I. Yakushin**}}

\vskip 1cm

\centerline{*  Dipto. di Scienze Fisiche, Universit\'a di Napoli,
80125 Napoli, Italy,}

\medskip

\centerline{** Dept. of Physics and Astronomy, Univ. of Alabama,
Tuscaloosa, Al 35487, U.S.A.}

\vskip 2cm

\vspace*{5mm}

\normalsize
\centerline{\bf ABSTRACT}
We write down the Poisson structure for a relativistic particle where
the Lorentz group does not act canonically, but instead
as a Poisson-Lie group.  In so doing we obtain the classical
limit of a particle moving on a noncommutative space possessing
$SL_q(2,C)$ invariance.  We show that if the standard
mass shell constraint is chosen for the Hamiltonian function,
then the particle interacts with the space-time.  We solve for the
trajectory and find that it originates and terminates at
singularities.

\vskip 2cm
\vspace*{5mm}
\scrollmode

\newpage
 \baselineskip=24pt

\section{Introduction}

A number of authors have investigated the effects of replacing
ordinary Lorentz invariance in quantum theory by $SL_q(2,C)$
invariance. [1-7]  The underlying space-time for such theories,
denoted by $M_q$,
is required to be noncommuting.  A consistent differential calculus
on $M_q$ has been obtained, and furthermore
 wave equations possessing $SL_q(2,C)$ invariance were found for
such a domain.\cite{Son}
  In addition, a corresponding
 $q-$deformed Poincare algebra was constructed,
along with their Hilbert space representations.\cite{PSW}

Here we shall
attempt to understand just what systems are being quantized in arriving
 at the above mentioned $SL_q(2,C)$ invariant quantum theories.  Thus
 we wish to examine the classical limit.
The classical analysis is considerably simpler
than its quantum counterpart because the classical space-time is ordinary
Minkowski space, and thus it is commuting,
and the classical
symmetries are associated with ordinary Lie groups, and not
quantum groups.    In addition, to check the
consistency of the algebra we essentially only need to verify the Jacobi
identity.

For us the relevant classical symmetry
 group is $SL(2,C)$, the covering group of
the Lorentz group.  Although it transforms the space-time coordinates
in the usual way, these transformations are not
canonical transformations.  Instead they are
Poisson Lie group transformations.[8-13]
A Poisson Lie group transformation is said to correspond to
 any group action on the classical observables which is a Poisson map.
In other words, a Poisson Lie group transformation preserves
the Poisson structure for the classical observables,
 and unlike a canonical transformation,
it carries its own Poisson structure.
This Poisson structure is then replaced by quantum group commutation
relations upon passing to the quantum theory.\cite{Tak}

In \cite{MSS} we illustrated
how the $SU(2)\times SU(2)$ chiral
symmetry of the isotropic rigid rotator could be implemented as a
 Lie-Poisson transformation.   In another example \cite{Zak},
the action of the Poincare group on two dimensional space-time
was realized by a Lie-Poisson transformation.
Here we shall show how to implement the Lorentz group on four dimensional
space-time in terms of a Lie-Poisson transformation.

We begin by examining $SL(2,C)$ Lie-Poisson
transformations on spinors in Sec. 2.
Quantum groups are known to have a covariant action on quantum spinors
(the coordinates of a quantum plane).\cite{Man}
  The analogous statement in classical mechanics is
 formulated in terms of Lie-Poisson transformations
 on spinors which satisfy certain quadratic Poisson bracket
relations.\cite{Skl}  $SL(2,C)$ group elements are then also
required to satisfy quadratic Poisson brackets.\cite{Tak}

In Sec. 3 we examine $SL(2,C)$ Lie-Poisson
transformations on vectors.  Like the spinors, the vectors are
required to satisfy quadratic Poisson brackets.
These brackets correspond to
 the classical limit of the reflection equations.\cite{AKR}
The results are extended in Sec. 4
to obtain a one parameter deformation of
the symplectic structure for a relativistic particle.
In order to satisfy Lie-Poisson invariance and the Jacobi identity,
we find it necessary to propose an infinite order polynomial
expression for the Poisson bracket between position and momentum.
We then examine the dynamics which follows from choosing the Hamiltonian
function
to be the standard mass shell constraint.   Unlike in the canonical
theory, this choice
does not yield a free particle since an interaction is present in
the Poisson structure.  However, as with the free particle,
 we are able to find ten conserved quantities,
the momenta and a deformation of the usual
angular momenta.  We  then use these conservation principles
to analytically solve for the particle trajectory and we find that it
originates and terminates at singularities.

Concluding remarks are made in Sec. 5.

\section{Lie-Poisson spinors}

Let $\psi=\pmatrix{u \cr v}$
be a Grassmann spinor associated with the $(0,{1\over2})$
of the Lorentz group.  Thus it transforms according to
\Eq{ \psi\rightarrow \psi' = g\psi \;,\quad g \in SL(2,C)\;.\label{trps}}
This transformation would be canonical if $u$ and $v$ were to
 satisfy trivial Poisson brackets.  We shall instead examine
the following quadratic relations:
\Eq{\pois{u}{v}=\pois{v}{u}=i\lambda uv\;,\quad
\pois{u}{u}=  \pois{v}{v}=0\;,\label{pbuv}  }
where $\lambda$ is a real parameter.  Now for $\lambda \ne 0,\;SL(2,C)$
does not act canonically as the Poisson structure defined by
(\ref{pbuv}) is not preserved under
 (\ref{trps}).   On the other hand, the Poisson brackets can be made
 to be invariant if we allow the matrix elements of $g$
to have certain nonvanishing brackets with each other (and vanishing
brackets with $\psi$).  Group multiplication is then said to be
a Poisson map and
such transformations are called Lie-Poisson.  If we write
 $$g=\pmatrix{a & b\cr c & d},\quad ad-bc=1\;,$$
then the Poisson bracket relations which are needed for this
purpose are
$$  \pois ab = - i \lambda ab\;,\quad
  \pois ac = - i \lambda ac\;,\quad
  \pois bd = - i \lambda bd\;, $$
\Eq{  \pois cd = - i \lambda cd\;,\quad
  \pois bc = 0\;,\quad
  \pois ad = -2i \lambda bc\;.\quad \label{pb3} }
{}From these relations it
 follows that  det $ g $ is a classical Casimir function, and hence
 (\ref{pb3}) is consistent with the constraint: det $g=1$.
Now from (\ref{pb3}) we have the result that the Poisson brackets
(\ref{pbuv}) are preserved under the action of $SL(2,C)$ and hence
such transformations are Lie-Poisson.

The Poisson brackets (\ref{pbuv}) and (\ref{pb3})
 can be expressed more compactly
if we introduce tensor product notation.  The Poisson structure
(\ref{pb3}) for $g$ can be written according to
\Eq { \{\gone,\gtwo\} = [\;r  \;,\;\gone \gtwo\;] \;,\label{lpb} }
where the indices $1$ and $2$ refer to two separate
vector spaces.  $r$, known as the classical r-matrix, acts
nontrivially on both vector spaces $1$ and $2$,
while $\gone=g\otimes \BI,\;\gtwo=\BI\otimes g ,$
where $\BI$ is the unit operator.
Written as a $4\times4$ matrix, $r$
takes the following form:
\Eq{
 r={{i\lambda}\over 2}\pmatrix{1 & & &  \cr &-1 & &  \cr &4 &-1& \cr
& & & 1\cr} \quad . }
It is known to satisfy the modified classical Yang-Baxter equation,
which here insures that the Jacobi identity is satisfied.
The Poisson structure defined by (\ref{lpb}) is preserved under
left and right group multiplication.  For this one defines an
$h\in SL(2,C)$ to satisfy the same relations as $g$,
\Eq { \{\hone,\htwo\} = [\;r  \;,\;\hone \htwo\;] \;, }
where $\hone=h\otimes \BI,\;\htwo=\BI\otimes h ,$
and in addition $ \{\hone,\gtwo\} = 0$.  Then for example under
left multiplication $g \rightarrow hg$, we get
\begin{eqnarray}
 \{\gone, \gtwo\}\rightarrow
 \{\hone \gone,\htwo \gtwo\}& = &[\;r  \;,\;\hone \htwo\;]\gone \gtwo +
\hone \htwo [\;r  \;,\;\gone \gtwo\;] \cr
& = &[\;r  \;,\;(\hone\gone)( \htwo\gtwo)\;]  \;.
\end{eqnarray}

The Poisson brackets (\ref{pbuv})
 can be rewritten according to:
\Eq{\pois{\psone}{\pstwo}={1\over 2} (r+r^\dagger)\psone\pstwo \;.
\label{rrda}}
where $\psone=\psi\otimes \BI,\;\pstwo=\BI\otimes \psi .$
Poisson brackets (\ref{rrda}) are symmetric and from them the Jacobi
identity trivially follows.  Using the fact that
$r-r^\dagger$ is an adjoint invariant, ie.
\Eq{[r-r^\dagger,\gone\gtwo] =0\;,\label{adin}}
 we can then show that
\Eq{\pois{\gone\psone}{\gtwo\pstwo}=
{1\over 2} (r+r^\dagger)(\gone\psone)(\gtwo\pstwo )\;,}
proving again that Lorentz transformations are Lie-Poisson.

In the above, $\lambda$ plays the role of a deformation parameter,
$\lambda \rightarrow 0$ corresponding to the canonical limit.
When $\lambda \rightarrow 0$, the right hand sides of (\ref{lpb}) and
 (\ref{rrda})
tends to zero as in the canonical theory, and also in that limit,
the Lie-Poisson transformations tend to canonical transformations.

\section{Lie-Poisson vectors}

We next apply the same procedure to the case of the
vector or $({1\over2},{1\over2})$ representation of $SL(2,C)$.
For this we can introduce the hermitean matrix
\Eq{  x= \pmatrix{-x_0-x_3 & -x_1+ix_2\cr -x_1-ix_2 & -x_0+x_3}\;,
 \label{xdef}}
$x_\mu$ denoting the space-time components of the vector.
$x$ transforms according to
\Eq{x\rightarrow x'= gxg^\dagger \;.\label{xtr} }
As before we will assume that the Poisson bracket of $g$ with itself
is given by (\ref{lpb}).  [We will also assume it has zero brackets
with $x$.]
Since the transformation (\ref{xtr}) involves $g^\dagger$ as well as $g$,
we will need to know the Poisson brackets of
$g^\dagger$ or $\bar g = (g^\dagger)^{-1}$.  For this we demand
that the Poisson structure for $g$ and $ g^\dagger$ is consistent
with complex conjugation, antisymmetry and the Jacobi identity.
All three of these conditions are met  for the following relations:
\begin{eqnarray}
 \{\gone,\gbtwo\}&=&[\;r  \;,\;\gone \gbtwo\;] \;,\label{ggb}  \\
 \{\gbone,\gtwo\}&=&[\;r^\dagger\;,\;\gbone \gtwo\;]\;,\label{gbg} \\
 \{\gbone,\gbtwo\}&=&[\;r^\dagger  \;,\;\gbone \gbtwo\;] \;,\label{gbgb}
\end{eqnarray}
Poisson brackets
(\ref{lpb}) and (\ref{ggb}-\ref{gbgb}) coincide with the classical limit
of the $SL_q(2,C)$ commutation relations given in refs.\cite{AKR}.
Using these relations we can see how the Poisson structure for $x$ is
transformed under (\ref{xtr}).  We find
\begin{eqnarray}
 \{\xone',\xtwo'\}&=&
 \{\gone \xone \gbone^{-1}, \gtwo \xtwo \gbtwo^{-1}\}  \cr
&=& \;\gone \gtwo \biggl(  \{\xone ,\xtwo \} -
r \xone \xtwo -\xone   \xtwo
r^\dagger  +\xtwo  r \xone  +\xone  r^\dagger \xtwo  \biggr) \gbone^{-1}
\gbtwo^{-1}  \cr
& &+\; r \xone'\xtwo'+\xone ' \xtwo'
r^\dagger  -\xtwo' r \xone' -\xone' r^\dagger \xtwo '
\end{eqnarray}
If for the Poisson bracket of $x$ with itself we now choose
 \Eq{ \pois{\xone}{\xtwo} =  r \xone\xtwo +\xone\xtwo
r^\dagger  -\xtwo r \xone -\xone r^\dagger \xtwo \;,\label{pbxx}}
then these brackets are preserved under Lorentz transformations
 and consequently (\ref{xtr}) is a Lie-Poisson transformation.
(More generally, we can add a term $\Sigma$ to (\ref{pbxx})
which Lorentz
transforms according to $\Sigma \rightarrow \gone  \gtwo \Sigma
\gbone^{-1} \gbtwo^{-1}$.  For simplicity, we shall not
consider this possibility here.)
It can be checked that
 the Poisson structure defined by (\ref{pbxx}) is consistent
with invariance under hermitean
conjugation, antisymmetry and the Jacobi identity.

In terms of space-time coordinates, (\ref{pbxx}) can be expressed by
\Eq{\pois{x_i}{x_j}=2\lambda \epsilon_{ijk} x_k (x_0+x_3)\;,\quad
\pois{x_i}{x_0}=0\;.}
Thus the time component $x_0$ is in the center of the algebra.
Also in the center is the spatial distance $\sqrt{x_ix_i}$
and consequently the
invariant space-time interval det $x$.  Analogous central
elements appear in the quantum theory, indicating that
simultaneous measurement of time, radial distance from some
origin, as well as, the projection along one spatial axis
is possible.

Using the fact that $r-r^\dagger$ is an adjoint invariant,
we can rewrite (\ref{pbxx}) according to
 \Eq{ \pois{\xone}{\xtwo} = (\mu r +\nu r^\dagger)
     \xone\xtwo +\xone\xtwo (\mu r^\dagger+\nu r)
 -\xtwo r \xone -\xone r^\dagger \xtwo \;,\label{pbxx2}}
where $\mu$ and $\nu$ are real and subject to the constraint
$\mu+\nu=1$.
Yet another way to write the Poisson structure is to introduce a new
hermitean matrix  $\tilde x$ which can be
defined in terms of $x$ by
\Eq{\tilde x= \sigma_2 x^T \sigma_2\;,}
 $T$ denoting transpose and $\sigma_2$ being the second Pauli matrix.
 It transforms according to
\Eq{   \tilde x \rightarrow \tilde x '=
 \bar g\tilde xg^{-1} \;.\label{trxt}}
This transformation is Lie-Poisson when $\tilde x$ satisfies
\Eq{\pois{\xtone}{\xttwo}  =   r \xtone\xttwo +\xtone\xttwo
r^\dagger  -\xttwo r^\dagger \xtone -\xtone r \xttwo \;,\label{pbxtxt}}
which is identical to (\ref{pbxx}).

As was the case with spinors,
$\lambda$ plays the role of a deformation parameter.
All Poisson brackets vanish when $\lambda \rightarrow 0$,
and the Lie-Poisson transformations tend to canonical transformations.

\section{Lie-Poisson particle dynamics}

We now let $x$ correspond to the space-time coordinates of a particle,
$x_i$ denoting the space coordinates.
The quantum analogue of (\ref{pbxx}), known as
the reflection equation\cite{AKR}, implies that the space coordinates
 do not commute
among themselves and therefore that the quantum mechanical
particle is nonlocalizable.

What about the particle momenta $p$?  We shall assume that,
like $x$, it is a Lorentz vector and thus
we know how to write down
its Poisson brackets with itself
so that Lorentz transformations are Lie-Poisson.
Here we find it more convenient to use
$\tilde p= \sigma_2 p^T \sigma_2\;,$ whose transformation is
 analogous to (\ref{trxt}).
Thus \Eq{\pois{\ptone}{\pttwo}  =   r \ptone\pttwo +\ptone\pttwo
r^\dagger  -\pttwo r^\dagger \ptone -\ptone r \pttwo \;.\label{pbptpt}}

More difficult to obtain are the
Poisson brackets between position and momenta.
The result should give some deformation of the canonical symplectic
structure.  Thus we want to recover
$\pois{x_\mu}{p_\nu} = \eta_{\mu\nu},\;\eta =diag (-1,1,1,1)$ in the
limit $\lambda \rightarrow 0$.
This is in addition to demanding that the bracket be preserved under
Lorentz transformations.
 It must also be antisymmetric, satisfy the Jacobi identity and be
 invariant under hermitean conjugation.

We first address the issue of Lorentz covariance.  Under a Lorentz
transformation
\begin{eqnarray}
 \{\xone',\pttwo'\}&=&
 \{\gone \xone \gbone^{-1}, \gbtwo \pttwo \gtwo^{-1}\}  \cr
&=& \;\gone \gbtwo \biggl(  \{\xone ,\pttwo \} -
r \xone \pttwo -\xone   \pttwo
r^\dagger +\pttwo  r \xone  +\xone  r^\dagger \pttwo  \biggr) \gbone^{-1}
\gtwo^{-1}  \cr
& &+\; r \xone'\pttwo'+\xone ' \pttwo'
r^\dagger  -\pttwo' r \xone' -\xone' r^\dagger \pttwo '   \;.
\end{eqnarray}
Lorentz covariance follows for
 \Eq{ \pois{\xone}{\pttwo} =  r \xone\pttwo +\xone\pttwo
r^\dagger  -\pttwo r \xone -\xone r^\dagger \pttwo + \Delta
\;,\label{pbxpt}}
provided that $\Delta$ transforms according to
\Eq{\Delta \rightarrow \Delta'=
 \gone \gbtwo \Delta\gbone^{-1} \gtwo^{-1} \;.
\label{Dt}}
The term $\Delta$ cannot be zero.  This is evident since
in order to recover canonical Poisson bracket relations
in the limit $\lambda \rightarrow 0$,
$\Delta$ must approach a constant matrix in that limit.  For $\Delta$ to
satisfy (\ref{Dt}), the constant matrix
can only be proportional to the permutation
matrix ${\cal P}$.  ${\cal P}$ is defined to
switch the vector spaces $1$ and $2$ upon
commutation.  Thus, eg. $\gone{\cal P}={\cal P}\gtwo$.
The normalization
\Eq{ \Delta\rightarrow -2 {\cal P}\;, \quad {\rm as} \quad
 \lambda\rightarrow 0 \;, \label{thec} }
gives the correct canonical limit.

In addition to satisfying (\ref{Dt}) and (\ref{thec}),
$\Delta$ must be hermitean.    To search for a consistent
solution when $\lambda \ne 0$ we choose the following ansatz:
 \Eq{ \Delta=-(\fone {\cal P} + {\cal P}\fone^\dagger)
 \;,\label{Dif}}  where $f$ is a $2\times 2$ matrix-valued function of
$x$ and $p$.  $f$ also in general depends on $\lambda$.
{}From (\ref{thec}),
it tends to the unit matrix $\BI$ when $ \lambda\rightarrow 0.$
{}From (\ref{Dt}), we have that $ f$ transforms according to
\Eq{  f \rightarrow  g      f g^{-1} \;,\label{flptr} }
under the action of the Lorentz group.
In general, $ f $ should be expressible as a polynomial
in $x \tilde p$, $p \tilde x$, $x \tilde x$ and $ p \tilde p$.
Presumably, it is also a polynomial in $\lambda$, with the
lowest order term equal to $\BI$.
In order to fix the higher
order terms, we now require that the brackets (\ref{pbxx}),
(\ref{pbptpt}) and (\ref{pbxpt}) satisfy the Jacobi identity.
After some work (and with the aid of algebraic manipulation packages)
we find that an infinite order polynomial expression for
$f$ is needed for the job.  The lowest order terms are given by:
\Eq{f=\BI + i\lambda x\tilde p -{{\lambda^2}\over 2} (x\tilde p)^2
 - {{\lambda^4}\over 8} (x \tilde p)^4 -
 {{\lambda^6}\over {16}} (x \tilde p)^6 - ...\quad.\label{fexpa}}
With $f$ so defined up to $6$ orders in $\lambda$, the Jacobi
identity is violated at the $8^{th}$ order in $\lambda$.  The
first five terms in the series expansion
suggests that $f$ has the following exact expression:
\Eq{
f= i\lambda x \tilde p +\sqrt{\BI - {\lambda^2} (x \tilde p)^2}  \;.}
To uniquely define the square root we demand that it goes to $\BI$
when $\lambda \rightarrow 1$.  Alternatively, we can write
 \Eq{f= \exp \{i\lambda
J \}\;,\quad {\rm where} \quad \sin \lambda J=\lambda x\tilde p\;.
\label{defJ}   }

Eqs. (\ref{pbxx}), (\ref{pbptpt}) and (\ref{pbxpt}) define a consistent
Poisson structure which is a deformation
of the canonical symplectic structure for a relativistic particle and
which is also Lie-Poisson with respect to Lorentz transformations.
We next remark on some of its consequences for relativistic particle
dynamics.

To proceed we must specify the Hamiltonian function.
A natural choice is the mass shell constraint
\Eq{ H=\alpha( \det p - m^2) =\alpha( \det \tilde p - m^2)
 \;,\label{fHam}}
where $\alpha$ is a Lagrange multiplier and $m$ denotes the mass.
As usual,
fixing $\alpha$ corresponds to selecting an evolution parameter $\tau$
for the particle trajectory.

 {\it Now because the Poisson structure for the particle is
not the canonical one, the Hamiltonian
(\ref{fHam}) does not describe a free
particle.}  To obtain the equations of motion for this system
we need the Poisson brackets of $\det \tilde p$ with $x$ and $\tilde p$.
 From (\ref{pbxx}), (\ref{pbptpt}) and (\ref{pbxpt}) we get
\begin{eqnarray}  \{x, \det \tilde p \}&=&-(f p + p f^\dagger)  \;,  \\
 \{\tilde p, \det \tilde p \}&=& 0 \;. \end{eqnarray}
We thus get
\begin{eqnarray} \dot x ={d\over {d\tau}} x&=& -\alpha ( e^{i\lambda
J} p +  p e^{-i \lambda J^\dagger} )\;,
 \label{xeom} \\ \dot p= {d\over {d\tau}} p&=& 0 \;. \end{eqnarray}

As in the canonical case the four components of
 momentum are conserved.  In addition to the conservation of
momentum, there is an analogue to the
 conservation of angular momentum.  To see this let us compute
the $\tau$ derivative of $J$ defined in (\ref{defJ}).  We first
note that
\Eq{\dot x \tilde p = -\alpha (f p + p f^\dagger) \tilde p
 =  -2\alpha m^2 \sqrt{\BI - {\lambda^2} (x \tilde p)^2}  \;, }
where we have used the mass shell constraint $\tilde p p=p \tilde p =
m^2 \BI$.  It follows  that $\dot x \tilde p$ commutes with
$x \tilde p$ and therefore we are justified in writing
\Eq{ \frac{d}{d\tau} J=\frac{1}{\lambda} \frac {d}{d\tau} \sin^{-1}
 (\lambda x\tilde p )  = \frac { \dot x \tilde p} {
\sqrt{\BI - {\lambda^2} (x \tilde p)^2} } =-2\alpha m^2 \BI \;.}
Thus only the trace of $J$ has a nonzero time derivative.  The
remaining components \Eq{j = J-\frac 1 2 \BI \;
{\rm Tr} J \label{angdef}}
are constants of the motion.  There are six real components in $j$,
and in the limit $\lambda \rightarrow 0$ they reduce to
$x\tilde p-\frac 1 2 \BI \; {\rm Tr} x\tilde p$ which correspond
to the canonical
angular momenta.  $j$ is therefore a deformation of the standard
angular momentum.

Using these constants of motion we can now integrate the equation
of motion (\ref{xeom}) for $x$,
\Eq{ \dot x =  \frac{ \dot\beta(\tau) }{4 m^2}
 ( e^{\frac i2 \lambda \beta(\tau)}k +
 e^{-\frac i2 \lambda\beta(\tau)} k^\dagger ) \;,
 \quad k=e^{i\lambda j} p \;,\quad \beta(\tau)={\rm Tr} J (\tau)
 \;,} and thus up to an overall constant hermitean matrix \Eq{  x(\tau)
 =  \frac 1{2 i \lambda m^2} ( e^{\frac i2 \lambda \beta(\tau)}k -
 e^{-\frac i2 \lambda\beta(\tau)} k^\dagger )+C\;,\label{xsol}}
where $C$ is a constant hermitean matrix.
Here we have used the result that $\beta$ is real which follows
from the definition (\ref{defJ}) of $J$ in terms of hermitean matrices
$x$ and $\tilde p$.   From (\ref{xdef}), the time coordinate $x_0$
is given by
\Eq {x_0 =- \frac 12 {\rm Tr} X= -\frac{|{\rm Tr} k|}
{2\lambda m^2}\sin \biggl(\frac {\lambda \beta(\tau)}2+\phi\biggr)
- \frac 12 {\rm Tr} C \;,}
where $\phi$= arg ${\rm Tr}k$.  It follows that the particle
has a finite lifetime equal to ${|{\rm Tr} k|}/(|\lambda| m^2)$.
For the example where $p$ points in the time-like direction and
$j$ is zero, then the lifetime is simply
 $2/(|\lambda | m)$.  [It is curious that if one compares this
classical result with the Heisenberg
uncertainty relation for virtual particles,
then one finds an upper bound for $|\lambda|$, $|\lambda|\le 2/\hbar$.]
If we now introduce a rescaled time $t$ defined by
\Eq{ t= \frac{\lambda m^2}{|{\rm Tr} k|} \biggl(x_0 + \frac12{\rm Tr}C
\biggr) +\frac12\;,\quad 0\le t\le 1\;,} then the solution
(\ref{xsol}) can be expressed by
\Eq{  x(t)  =   U\; t^{1/2}\sqrt{1 - t } \;+\;Vt+x(0) \;,} where
$U$ and $V$ are constant
hermitean matrices ($U$ also being traceless) such that
\Eq{U-iV=-\frac2{\lambda m^2}\;ke^{-i\phi} \;.}

The trajectory is a superposition of
a semicircular path plus a line segment.  It
originates and terminates at singularities.  This because
the velocity  $\frac{dx}{dt}$  is singular when
$t=0,1$.   It  equals $V$ when $t=\frac12$.
Since the lifetime is finite (and the
speed can exceed the speed of light at some points on the trajectory),
the solutions at best can be interpreted as describing
virtual particles.  [On the other hand,
we can avoid problems associated with superluminal speeds if
we replace $x_i$ with new spatial coordinates $X_i$, where
$$X_i =x_i - U_i \; t^{1/2}\sqrt{1 - t } \;\;.$$
 In terms of $X_i$
 the particle trajectory is a straight line segment.
However for such coordinates Lorentz covariance is lost and
the Poisson structure becomes more complicated.]

\section{Concluding Remarks}

Here we summarize our results and indicate future avenues of
research.

We have found a Poisson structure for the relativistic particle
where the Lorentz group is implemented as a Lie-Poisson
transformation.  It involves an infinite order polynomial function
$f$ in the variables $x$ and $p$.  We have obtained the lowest order
terms in the expansion for $f$ in (\ref{fexpa}) and surmised that the
exact expression is given by (\ref{defJ}).  It would be desirable to
verify that this exact expression is consistent with the Jacobi
identity.

As in the canonical theory of a relativistic free particle,
our dynamical system has ten conserved quantities $p$ and $j$. When
$\lambda\rightarrow 0$, these quantities reduce to
the canonical momenta and angular momenta.   For $\lambda \ne 0$,
the Poisson brackets of $p$ and $j$
 give a deformation of
the Poincare algebra.  This deformed Poincare algebra has the
property that one of its Casimirs is the ``mass"-squared, det $p$.
In order to understand what spin means in this theory it would
be useful to know what other Casimirs (if any) exist.
If we replace our Hamiltonian (\ref{fHam})
by a more general function of the Casimirs, we can
obtain a large family of dynamical systems all having $p$ and $j$ as
conserved quantities.  Perhaps one such system will give a
noncanonical description of the standard relativistic
free particle.  One may then perform a deformation quantization
of that system thereby giving an alternative quantization of the
free particle where  $SL_q(2,C)$ corresponds to
 the space-time symmetry.

With regards to the deformed angular momenta $j$, it would be
of interest to know whether or not $j$ can be interpreted
as the generator of the Lie-Poisson Lorentz symmetry.
In the study of the Lie-Poisson formulation of the
rigid rotator in \cite{MSS} the generators of the Lie-Poisson
chiral symmetry were found to have
unusual properties (as opposed to canonical symmetry generators).
For example, they took values in a group dual to the symmetry
group.  Analogous novel features are
anticipated here.  It is clear for instance
that infinitesimal Lorentz
transformations are not obtained (as in the canonical theory)
simply by taking Poisson
brackets with the conserved angular momenta
$j$.   Similarly, translations are not obtained
by taking Poisson brackets with $p$.  In this regard, it would be of
interest to check whether or not translations can be implemented in
a Lie-Poisson fashion.  If so, the entire Poincare
group may be a Lie-Poisson symmetry of this theory, thus yielding
 a four dimensional analogue of the two dimensional
Poincare invariant theory worked out in \cite{Zak}.

We plan to address the above issues along with quantization
in a future article.
\bigskip
\bigskip

{\parindent 0cm{\bf Acknowledgements:}}
We are very grateful to G. Marmo for useful discussions and for informing
us of ref. \cite{Zak}.
 A. S. would like to thank the members of
the theory group of the University of Naples for their
 warm  hospitality while this work was being initiated.
 A. S. was supported in part
by the Department of Energy, USA, under contract number
DE-FG05-84ER40141.

\end{document}